\documentclass[a4paper,11pt]{article}
\usepackage{graphicx}
\usepackage{amsmath}
\usepackage{amssymb}
\usepackage{float}

\newcommand{\beq}{\begin{equation}}
\newcommand{\eeq}{\end{equation}}

\newcommand{\ba}{\begin{array}}
\newcommand{\ea}{\end{array}}

\newcommand{\bea}{\begin{eqnarray}}
\newcommand{\eea}{\end{eqnarray}}

\newcommand{\bc}{\begin{center}}
\newcommand{\ec}{\end{center}}

\newcommand{\bt}{\begin{tabular}}
\newcommand{\et}{\end{tabular}}

\newcommand{\bi}{\begin{itemize}}
\newcommand{\ei}{\end{itemize}}

\newcommand{\bd}{\begin{description}}
\newcommand{\ed}{\end{description}}

\newcommand{\bp}{\begin{pmatrix}}
\newcommand{\ep}{\end{pmatrix}}

\newcommand{\pd}{\partial}

\newcommand{\gnorm}[1]{\left|\left| #1\right|\right|}

\title{The Dynamic-Mode Decomposition and Optimal Prediction}
\author{Christopher W. Curtis \quad\quad Daniel Jay Alford-Lago}
\date{}
\begin{document}
\maketitle

\begin{abstract}
The Dynamic-Mode Decomposition (DMD) is a well established data-driven method of finding temporally evolving linear-mode decompositions of nonlinear time series. Traditionally, this method presumes that all relevant dimensions are sampled through measurement. To address dynamical systems in which the data may be incomplete or represent only partial observation of a more complex system, we extend the DMD algorithm by including a Mori-Zwanzig Decomposition to derive memory kernels that capture the averaged dynamics of the unresolved variables as projected onto the resolved dimensions. From this, we then derive what we call the Memory-Dependent Dynamic Mode Decomposition (MDDMD). Through numerical examples, the MDDMD method is shown to produce reasonable approximations of the ensemble-averaged dynamics of the full system given a single time series measurement of the resolved variables. 
\end{abstract}

%%%%%%%%%%%%%%%%%%%%%%%%%%%%%%%%%%%%%%%%%%%%%%%%%%%%%%%%%%%%%%%%%%%%%%%%%%%%%%%%
% Introduction
%%%%%%%%%%%%%%%%%%%%%%%%%%%%%%%%%%%%%%%%%%%%%%%%%%%%%%%%%%%%%%%%%%%%%%%%%%%%%%%%
\section{Introduction}

In recent years, beginning with the now seminal papers \cite{mezic1,schmid}, a large body of literature has emerged around the modal decomposition technique known as the Dynamic-Mode Decomposition (DMD); see \cite{kutz} for a comprehensive review of the literature, and the recent reviews in \cite{taira} and \cite{clainche} for the most recent summaries and extensions.  The primary advantage of DMD is that it is a completely model free data processing tool, allowing for the ready determination of results from arbitrarily complicated data sets.  Moreover, in contrast to other methods such as Principle Component Analysis, the DMD also provides a ready means for generating models from measured data alone, thereby allowing for a means of sidestepping model development in cases where it would either be difficult or even impossible to do so.  However, this advantage comes at the cost of introducing an infinite dimensional framework by way of determining the spectral properties of the affiliated Koopman operator \cite{koopman}.  Authors have, through using functional analytic ideas, been able to exploit this infinite dimensional structure to enhance the performance of the DMD \cite{williams,williams2}.  However, it has been shown that these enhancements are difficult to apply in a systematic way \cite{kutz2}, and thus approaches merging the DMD with nonlinear optimization are currently being explored to improve the accuracy and usability of the method \cite{bertozzi,kutz3}.

In a very different direction, building on the now seminal work from the statistical mechanics community in \cite{mori} and \cite{zwanzig}, in \cite{chorin,chorin2,chorin3}, by way of what was called a Mori--Zwanzig Decomposition (MZD), the Koopman operator formalism was ultimately shown to be an effective means for accounting for missing information in measurements related to dynamical systems.  In these works, formulas for {\it memory kernels} were derived whereby the impact of unsampled dimensions in a dynamical system could be accounted for in a statistically consistent way, thereby helping to better predict the correct averaged results one would expect if one had an essentially unlimited number of measurements.  Given that this is never the case in practical situations, the MZD approach then provides an especially attractive means for approximating averaged dynamics in complex dynamical systems.  

Other authors have used this MZD approach to help generate more refined coarse-scale molecular dynamics simulations \cite{karniadakis}, close numerical schemes with respect to sub-grid scale variables \cite{parish2017unified}, and help provide long-time accurate reduced order models for nonlinear dynamical systems \cite{stinis}.  Motivated by this, in this work we merge the Koopman operator based MZD method with the standard DMD algorithm to introduce memory dependence into the DMD algorithm.  We term this new algorithm the Memory Dependent DMD (MDDMD).  Using this, we are able to derive a model-independent algorithm which allows from a single time series of data, representing only a partial measurement of a dynamical system, the accurate estimation of the reduced dynamics found by averaging over the missing or unmeasured dimensions.  Our results are supported with numerical simulation, which show our method to be quite effective.  Thus, we have derived a means for predicting averaged dynamics from one partial time series without making recourse to any particular dynamical system or model.  Likewise, as with the standard DMD, we are able to generate a corresponding model which could be used for prediction beyond the given time series.  The utility of such a feature will be explored in future work.  

The rest of the paper is structured in the following way.  The remainder of this section is devoted to an introduction to the Koopman operator, the standard DMD method, and the MZD method.  Section 2 presents a derivation of our MDDMD algorithm.  Section 3 presents numerical examples and their results.  Finally, Section 4 presents our conclusions and discussions of future work.  An Appendix summarizes technical results used throughout the paper.  

%%%%%%%%%%%%%%%%%%%%%%%%%%%%%%%%%%%%%%%%%%%%%%%%%%%%%%%%%%%%%%%%%%%%%%%%%%%%%%%%
\subsection{Koopman Operators and the Dynamic-Mode Decomposition}
We study nonlinear dynamical systems of the generic form
\[
\frac{d}{dt}{\bf y} = f({\bf y}), ~ {\bf y}(0) = {\bf x} \in \mathbb{R}^{n},
\]
where we suppose $N_{s}\gg1$.  We denote the affiliated flow of this equation as ${\bf y}(t) = \varphi(t;{\bf x})$.  We define the associated Hilbert space of {\it observables}, say $L_{2}\left(\mathbb{R}^{N_{s}},\mathbb{R},\mu\right)$, so that $g \in L_{2}\left( \mathbb{R}^{N_{s}},\mathbb{R},\mu\right)$ if
\[
\int_{\mathbb{R}^{N_{s}}} \left|g({\bf x})\right|^{2} d\mu\left({\bf x}\right) < \infty,
\]
where $\mu$ is some appropriately chosen measure.  A great deal of insight can be gained from looking at the affiliated linear representation of the problem given via the infinite-dimensional Koopman operator $\mathcal{K}^{t}$, where 
\[
\mathcal{K}^{t}:L_{2}\left( \mathbb{R}^{N_{s}},\mathbb{R},\mu\right)\rightarrow L_{2}\left( \mathbb{R}^{N_{s}},\mathbb{R},\mu\right) 
\]
so that if $g \in L_{2}\left( \mathbb{R}^{N_{s}},\mathbb{R},\mu\right)$, then 
\[
\mathcal{K}^{t}g({\bf x}) =  g\left(\varphi(t;{\bf x})\right).
\]

The power in this approach is that by moving to a linear-operator framework, the affiliated dynamics of the nonlinear system as measured through observables is captured by the eigenvalues of $\mathcal{K}^{t}$.  Thus, assuming for the moment that $\mathcal{K}^{t}$ has only discrete spectra, then if we can find a basis of $L_{2}\left( \mathbb{R}^{N_{s}},\mathbb{R},\mu\right)$ via the Koopman modes $h_{j}$ where
\[
\mathcal{K}^{t}h_{j} = e^{t\lambda_{j}}h_{j},
\]
then for any other observable $g$, we should have 
\[
g = \sum_{j=1}^{\infty}c_{j}h_{j}, ~ \mathcal{K}^{t}g = \sum_{j=1}^{\infty}e^{\lambda_{j}t}c_{j}h_{j}.
\]
However, as one would imagine, determining the modes of the Koopman operator is in general impossible in an analytic way.  The Dynamic-Mode Decomposition (DMD) method \cite{schmid,mezic1,williams,kutz} has been developed which allows for the approximation of a finite number of the Koopman modes.  

To wit, if we sample the flow $\varphi(t;{\bf x})$ at a discrete set of times $t_{j}=j\delta t$, where $j=0,\cdots, N_{T}$, thereby generating the data set ${\bf y}_{j} = \varphi(t_{j},{\bf x})$, and if we select some set of observables ${\bf g} = \{g_{l}\}_{l=1}^{M}$, then the DMD method approximates $\mathcal{K}^{\delta t}$ by computing the spectra of the finite-dimensional operator $\tilde{\mathcal{K}}_{a}$, which solves the finite-dimensional least-squares problem in the Frobenius norm 
\begin{equation}
\tilde{\mathcal{K}}_{a} = \mbox{arg min}_{A}\gnorm{G_{+} - A G_{-}}_{F}
\label{dmdeq}
\end{equation}
where $G_{-}$ and $G_{+}$ are the $M\times N_{T}$ matrices of the past and future states respectively, i.e.
\[
G_{-}=\left\{{\bf g}\left({\bf y}_{0}\right) \cdots {\bf g}\left({\bf y}_{N_{T}-1}\right) \right\}, ~ G_{+} = \left\{{\bf g}\left({\bf y}_{1}\right) \cdots {\bf g}\left({\bf y}_{N_{T}}\right) \right\}
\]
Note, the Frobenius norm $\gnorm{\cdot}_{F}$ of a matrix $A$ is given by 
\[
\gnorm{A}_{F}^{2} = \mbox{tr}\left(A^{\dagger}A \right),
\]
where $\mbox{tr}$ is the trace of a matrix and the $^{\dagger}$ denotes the Hermitian conjugate of a matrix.

Using a Singular-Value Decomposition (SVD) on $G_{-}$, so that $G_{-} = U\Sigma V^{\dagger}$, where $\Sigma$ is diagonal with strictly non-zero entries reflecting the rank of $G_{-}$, and where $U$ and $V$ are reductions of unitary matrices such that  
\[
U^{\dagger}U = V^{\dagger}V = I_{N_{r}},  
\]
where $N_{r}$ is the rank of $G_{-}$, and $I_{N_{r}}$ is the $N_{r}\times N_{r}$ identity matrix.  Note, the $^{\dagger}$ denotes the Hermitian conjugate of a matrix. Thus we can solve the minimization problem by letting
\[
\tilde{K}_{a} = G_{+}G^{-P}_{-} = G_{+}V\Sigma^{-1}U^{\dagger},
\]
where $G_{-}^{-P}$ denotes the Moore--Penrose pseudoinverse.  

%%%%%%%%%%%%%%%%%%%%%%%%%%%%%%%%%%%%%%%%%%%%%%%%%%%%%%%%%%%%%%%%%%%%%%%%%%%%%%%%%%%
\subsection{The Liouville Equation and the Mori--Zwanzig Decomposition}
Treating the nonlinear dynamical system as a vector field, we can form the affiliated Liouville equation which is an initial value problem of the form 
\[
u_{t} = \mathcal{L} u, ~ u({\bf x},0) = g({\bf x}),
\]
where
\[
\mathcal{L} = \sum_{j=1}^{N}f_{j}({\bf x})\pd_{x_{j}}
\]
so that 
\[
u({\bf x},t) = e^{\mathcal{L}t}g({\bf x}) = g\left(\varphi\left(t;{\bf x}\right)\right) = \mathcal{K}^{t}g({\bf x}).
\]
Thus we see that the Koopman operator is a stand in for the affiliated semi-group $e^{\mathcal{L}t}$.  To wit, if
\[
\dot{\varphi}_{j} = f_{j}(\varphi(t,{\bf x})), 
\]
then letting $g({\bf x}) = x_{j}$, we have 
\[
\pd_{t}\left(e^{\mathcal{L}t}x_{j} \right)  = e^{\mathcal{L}t}\mathcal{L}x_{j} = e^{\mathcal{L}t}f_{j}\left({\bf x}\right).
\]
We now suppose that there is an orthogonal projection $\mathbb{P}$ acting on $L_{2}\left( \mathbb{R}^{N_{s}},\mathbb{R},\mu\right)$ with corresponding complement $\mathbb{Q} = I - \mathbb{P}$.  Through an approach first derived in \cite{mori} and \cite{zwanzig} and generalized in \cite{chorin}, we are able to rewrite the given dynamical system or equivalent Liouville equation in the form 
\begin{equation}
\pd_{t}\left(\mathbb{P}e^{\mathcal{L}t} g \right) = \mathbb{P} e^{\mathcal{L}t}\mathbb{P}\mathcal{L}g + \int_{0}^{t}\mathbb{P}e^{\mathcal{L}(t-s)}K(s;g) ~ds
\label{mzeq1}
\end{equation}
and
\begin{equation}
K(t;g) + \int_{0}^{t}\mathbb{P}e^{\mathcal{L}(t-s)}\mathcal{L}K(s;g) ~ ds= \mathbb{P}e^{\mathcal{L}t}\tilde{F}(0;g),
\label{mzeq2}
\end{equation}
where $K(t;g)$ is called the {\it memory kernel}, given by 
\[
K(t;g) = \mathbb{P}\mathcal{L}F(t;g),
\]
where the {\it noise} $F(t;g)$ is given by the dynamics orthogonal to $\mathbb{P}$ so that 
\[
F(t;g) =  e^{\mathbb{Q}\mathcal{L}t}\mathbb{Q}\mathcal{L}g.
\]
Note, Equations \eqref{mzeq1} and \eqref{mzeq2}, hereafter known as the Mori--Zwanzig decomposition (MZD) are not approximations, instead being an equivalent reformulation of the Liouville initial value problem.  The focus of this paper is on discretizing the MZD in such a way as to only depend on knowledge of the observations $g(\varphi_{t}({\bf x}))$ with no recourse made to $\mathcal{L}$ and thus an explicit model dynamical system $\dot{{\bf y}} = f\left({\bf y}\right)$.  

%%%%%%%%%%%%%%%%%%%%%%%%%%%%%%%%%%%%%%%%%%%%%%%%%%%%%%%%%%%%%%%%%%%%%%%%%%%%%%%%
% Memory Dependent DMD: the DMD and Mori--Zwanzig Decomposition
%%%%%%%%%%%%%%%%%%%%%%%%%%%%%%%%%%%%%%%%%%%%%%%%%%%%%%%%%%%%%%%%%%%%%%%%%%%%%%%%
\section{Memory Dependent DMD: the DMD and Mori--Zwanzig Decomposition}
We now suppose that a given observable $g({\bf x})$ is dependent on two different, orthogonal sets of coordinates so that $g({\bf x}) = g(\hat{{\bf x}},\tilde{{\bf x}})$, where we imagine that $\hat{{\bf x}}$ represents those coordinates which are accessible to measurement while $\tilde{{\bf x}}$ are not readily observable.  Given an appropriate probability measure, say $\tilde{p}$, over the unmeasurable variables $\tilde{{\bf x}}$, then for any observable $g({\bf x},t)$ we have a projection $\mathbb{P}$ defined through conditional expectations, i.e.  
\[
\mathbb{P}g(\hat{{\bf x}}) = \mathbb{E}\left[g({\bf x})\left|\hat{{\bf x}}\right.\right] = \int g(\hat{{\bf x}},\tilde{{\bf x}}) \tilde{p}(d\tilde{{\bf x}}).
\]
Thus, we suppose that our `data' stream is given by $\hat{g}_{n}\in \mathbb{C}^{N_{d}}$ where
\[
\hat{g}_{n}(\hat{{\bf x}}) = \mathbb{P}e^{\mathcal{L}t_{n}}g, ~ t_{n} = n\delta t, 
\]
i.e. all we can measure are the time dynamics of a process with regards to the observable dimensions $\hat{{\bf x}}$.   Likewise, we treat $\mathcal{L}$ and $e^{t\mathcal{L}}$ as terms to discover from the data.  Following the DMD paradigm, we define 
\[
\tilde{T}(t) = \mathbb{P}e^{\mathcal{L}t}.
\]
Ultimately, we wish to find a finite-dimensional diagonalization of $\tilde{T}(\delta t)$ so that if 
\[
\tilde{T}(\delta t) \approx Ve^{\delta t \Lambda} V^{-1}, 
\]
then
\[
\tilde{T}(t) \approx V e^{t \Lambda} V^{-1}.
\]
Thus, we have the consistent approximations
\[
\mathbb{P}\mathcal{L} \approx V \Lambda  V^{-1}, ~ \mathbb{P} e^{\mathcal{L}t_{n}}\mathbb{P}\mathcal{L}g \approx \mathbb{P}\mathcal{L} \hat{g}_{n}, ~ \tilde{F}(0;g) \approx K(0;g).
\]
In order to get a closed form approximation, in Equation \eqref{mzeq2}, we further suppose that 
\begin{equation}
\int_{0}^{t}\mathbb{P}e^{\mathcal{L}(t-s)}\mathbb{Q}\mathcal{L}K ds \approx 0,
\label{closeeq}
\end{equation}
and we use a first-order finite-difference approximation to the derivative.  Note, we see that the approximation in Equation \eqref{closeeq} is internally consistent with the assumptions we have made by seeing that 
\[
\mathbb{P}\mathcal{L}^{j} \approx V\Lambda^{j}V^{-1}.  
\]
Thus, using the expansion
\[
\mathbb{P}e^{\mathcal{L}(t-s)}\mathbb{Q}\mathcal{L} = \sum_{j=1}^{\infty}\frac{(t-s)^{j}}{j!}\left(\mathbb{P}\mathcal{L}^{j+1}-\mathbb{P}\mathcal{L}^{j}\mathbb{P}\mathcal{L} \right),
\]
we readily get Equation \eqref{closeeq}.  Note, this does not show that the DMD is somehow inexplicably accurate.  Instead, it shows that the DMD assumptions we have made prevent us from directly determining any aspect of the orthogonal dynamics associated with the projection $\mathbb{Q}$.  

Using the Trapezoid Method for the remaining integrals and letting 
\[
\tilde{g}_{n} = V^{-1}\hat{g}_{n}, ~ \tilde{K}_{n} = V^{-1}K(t_{n};g),
\]
we get that 
\begin{align*}
\tilde{g}_{n+1} = & \left(I + \delta t \Lambda\right)\tilde{g}_{n}+\frac{(\delta t)^{2}}{2}\left(\tilde{K}_{n} + e^{n\delta t \Lambda}\tilde{K}_{0}+2 \sum_{l=1}^{n-1}e^{(n-l)\delta t \Lambda}\tilde{K}_{l}\right),\\
\tilde{K}_{n} = & \left(I + \frac{\delta t}{2}\Lambda\right)^{-1}\left( e^{n \delta t \Lambda}\left(I - \frac{\delta t}{2}\Lambda\right) \tilde{K}_{0} - \delta t\sum_{l=1}^{n-1}e^{(n-l)\delta t \Lambda}\tilde{K}_{l} \right).
\end{align*}
From the discretized Volterra equation, we find the recursive sequence
\[
\tilde{K}_{n} = e^{\delta t \Lambda}M(\Lambda)\tilde{K}_{n-1}, ~ n\geq 1, 
\]
where
\[
M(\Lambda) = I - \delta t\left(I + \frac{\delta t}{2}\Lambda\right)^{-1}.
\]
We then get the update formula 
\begin{align*}
\tilde{g}_{n+1} = & \left(I + \delta t \Lambda\right)\tilde{g}_{n} + \frac{(\delta t)^{2}}{2}e^{n\delta t \Lambda}\left(M^{n}-I\right)\left(I + 2 (M-I)^{-1} \right)\tilde{K}_{0},\\
= & \left(I + \delta t \Lambda\right)\tilde{g}_{n} - \delta t e^{n\delta t \Lambda}\left(M^{n}-I\right)\left( 1 - \frac{\delta t }{2} + \frac{\delta t }{2}\Lambda\right)\tilde{K}_{0}.
\end{align*}
Letting 
\[
\tilde{\Lambda} = I + \delta t \Lambda, 
\]
and $A = V\tilde{\Lambda}V^{-1}$, so that
\[
\mathbb{P}\mathcal{L}\approx \frac{1}{\delta t}(A-I), 
\]
by using the spectral representation theorem for matrices, we get the extended, memory-dependent optimization problem 
\begin{equation}
\tilde{\mathcal{K}}_{a} = \mbox{arg min}_{A}\gnorm{G_{+} - A G_{-} +  \delta t  \tilde{M}(A;K_{0})}^{2}_{F},
\label{exdmdeq}
\end{equation}
where we have 
\[
\tilde{M}(A;K_{0}) = \left(1 - \frac{\delta t}{2}  + \frac{1}{2}(A-I) \right)\left\{0 ~ f_{1}(A;K_{0}) \cdots f_{n}(A;K_{0})\right\},
\]
with
\[
f_{j}(A;K_{0}) = e^{j(A-I)}\left(M^{j}(A)-I\right)K_{0},
\]
where
\[
M(A) = I - \delta t \left(I + \frac{(A-I)}{2} \right)^{-1}.
\]

This in turn gives us the nonlinear optimization problem for finding $\tilde{\mathcal{K}}_{a}$ so that 
\[
\tilde{\mathcal{K}}_{a}G_{-}G_{-}^{T} - G_{+}G_{-}^{T} +  \delta t  \left(R^{T}(\tilde{\mathcal{K}}_{a})D\tilde{M}\left(\tilde{\mathcal{K}}_{a};K_{0}\right)-G_{-}\tilde{M}^{T}\left(\tilde{\mathcal{K}}_{a};K_{0}\right)\right)+ \cdots= 0, 
\]
where $D\tilde{M}\left(\tilde{\mathcal{K}}_{a};K_{0}\right)$ is found by computing 
\[
\mbox{tr}\left(R^{T}(\tilde{\mathcal{K}}_{a})D\tilde{M} ~W \right) = \lim_{\epsilon\rightarrow 0} \mbox{tr}\left(R^{T}(\tilde{\mathcal{K}}_{a})\frac{\left(\tilde{M}\left(\tilde{\mathcal{K}}_{a}+\epsilon W;K_{0}\right)-\tilde{M}\left(\tilde{\mathcal{K}}_{a};K_{0}\right)\right)}{\epsilon} \right),
\]
where
\[
R(\tilde{\mathcal{K}}_{a}) = G_{+}-\tilde{\mathcal{K}}_{a}G_{-}.
\]
Expanding $\tilde{\mathcal{K}}_{a} = \tilde{\mathcal{K}}_{a,0} + \delta t  \tilde{\mathcal{K}}_{a,1} + \cdots$, we readily get that 
\[
\tilde{\mathcal{K}}_{a,0} = G_{+}G_{-}^{-P},
\]
and
\[
\tilde{K}_{a,1}(K_{0}) =  -\left(R^{T}(\tilde{\mathcal{K}}_{a})D\tilde{M}\left(\tilde{\mathcal{K}}_{a,0};K_{0}\right) -G_{-} \tilde{M}^{T}\left(\tilde{\mathcal{K}}_{a,0};K_{0}\right)\right)\left( G_{-}G_{-}^{T}\right)^{-1}.
\]
Note, see the Appendix for details on how to approximate $D\tilde{M}\left(\tilde{\mathcal{K}}_{a,0};K_{0}\right)$ for computational purposes.  

We now follow the standard DMD procedure and recover the traditional DMD decomposition by diagonalizing $\tilde{K}_{a,0}$ so that 
\[
\tilde{\mathcal{K}}_{a,0} = V_{0}\tilde{\Lambda}_{0}V^{-1}_{0}.
\]
Pursuing a perturbative result, we then look at the eigenvalue problem 
\[
\left(\tilde{\mathcal{K}}_{a,0} + \delta t \tilde{K}_{a,1}(K_{0}) \right) \left( {\bf v}_{0,j} + \delta t {\bf v}_{1,j}(K_{0}) \right)  \approx \left(\tilde{\lambda}_{0,j} + \delta t \tilde{\lambda}_{1,j}(K_{0}) \right)\left({\bf v}_{0,j} + \delta t{\bf v}_{1,j}\right),
\]
so that we find, assuming for now the simplicity of the unperturbed eigenvalue $\lambda_{0,j}$, that 
\begin{equation}
\tilde{\lambda}_{1,j}(K_{0}) = \frac{\left<\tilde{K}_{a,1}(K_{0}){\bf v}_{0,j},{\bf v}^{(a)}_{0,j}\right>}{\left<{\bf v}_{0,j},{\bf v}^{(a)}_{0,j}\right>},  ~ {\bf v}^{(a)}_{0,j} \in \mathcal{N}\left\{\left(\tilde{\mathcal{K}}_{a,0}-\lambda_{0,j}\right)^{\dagger} \right\}, 
\label{perteqs}
\end{equation}
where $\mathcal{N}$ denotes the null-space of a matrix, and $\left<\cdot,\cdot\right>$ is the canonical inner product on $\mathbb{C}^{N_{d}}$.   This likewise allows for one to find the perturbation to the eigenvector using the formula
\[
{\bf v}_{1,j}(K_{0}) = \left(\tilde{\mathcal{K}}_{a,0}-\tilde{\lambda}_{j,0} \right)^{-P}\left(\tilde{\lambda}_{1,j}(K_{0}) -  \tilde{K}_{a,1}(K_{0})\right){\bf v}_{0,j}.
\]

Thus, treating $K_{0}$ as a random vector, we see that we have now generated an ensemble of approximations to the projected semi-group $\tilde{T}(\delta t)=\mathbb{P}e^{\mathcal{L}\delta t}$.  Averaging over these then gives us the approximation 
\[
\tilde{T}(t) g \approx \sum_{l=1}^{N_{r}}\left< \left({\bf v}_{l,0} + \delta t {\bf v}_{l,1}(K_{0}) \right) \mbox{exp}\left(\frac{t}{\delta t}\left(\tilde{\lambda}_{0,l}-1+\delta t \tilde{\lambda}_{1,l}(K_{0}) \right)\right)\right>_{K_{0}} b_{l},
\]
where again $N_{r}$ is the rank of $\tilde{K}_{a,0}$.  Decomposing the stochastic terms into their means and fluctuations so that 
\[
{\bf v}_{l,1}(K_{0})  =  \bar{{\bf v}}_{l,1} + {\bf v}^{(f)}_{l,1}(K_{0}), 
\]
and
\[
\mbox{exp}\left(\delta t \tilde{\lambda}_{1,l}(K_{0}) \right) =\bar{\lambda}_{1,l} + \tilde{\lambda}^{(f)}_{1,l}(K_{0}),
\]
we then get up to $\mathcal{O}(\delta t)$ that 
\begin{equation}
\tilde{T}(t) g \approx \sum_{l=1}^{N_{r}} \left( {\bf v}_{l,0} + \delta t \bar{{\bf v}}_{l,1} \right) \mbox{exp}\left(\frac{t}{\delta t}\left(\tilde{\lambda}_{0,l}-1+\log \bar{\lambda}_{1,l}\right)\right)  b_{l}.
\label{mddmdeq}
\end{equation}
This then gives us a formula for finding what we call the memory dependent DMD (MDDMD).  We can summarize our algorithm in the following way.

\begin{enumerate}
\item Given data stream $G_{-}$ and $G_{+}$, compute and diagonalize 
\[
\mathcal{K}_{a,0} = G_{+}G^{-P}_{-} = V_{0}\Lambda_{0}V^{-1}_{0}.  
\]
This is the standard DMD algorithm.  
\item From a given distribution, generate an ensemble of random vectors $K_{0}$.  Use these to compute $\tilde{\lambda}_{j}(K_{0})$ and ${\bf v}^{(a)}_{0,j}$ as in Equation \eqref{perteqs}.
\item Compute the averages of these perturbation parameters and use them to construct the MDDMD approximation given by Equation \eqref{mddmdeq}.
\end{enumerate}

%%%%%%%%%%%%%%%%%%%%%%%%%%%%%%%%%%%%%%%%%%%%%%%%%%%%%%%%%%%%%%%%%%%%%%%%%%%%%%%%%%%
\section{Examples and Results}
\subsection{A Hamiltonian System} \label{example_1}
Following the choice of example in \cite{chorin}, we study the four-dimensional system generated by the Hamiltonian, $H({\bf y})$, where 
\[
H({\bf y}) = \frac{1}{2}\left( y_{1}^{2}+y_{2}^{2}+y_{3}^{2}+y_{4}^{2}+ y_{1}^{2}y_{3}^{2}\right),
\]
with conjugate pairs $(y_{1}(t),y_{2}(t))$ and $(y_{3}(t),y_{4}(t))$.  This system describes two nonlinearly coupled oscillators
\begin{align*}
\dot{y}_1 &= y_2, \\
\dot{y}_2 &= - y_1 \left(1 + y_3^2 \right), \\
\dot{y}_3 &= y_4, \\
\dot{y}_4 &= -y_3 \left(1 + y_1^2 \right).
\end{align*}
We treat the initial conditions $y_{1}(0)=x_{1}$ and $y_{2}(0)=x_{2}$ as known while the remaining initial conditions are considered unknown, and thus $\hat{{\bf x}}=(x_{1},x_{2})$ and $\tilde{{\bf x}}=(x_{3},x_{4})$.  Fixing $\hat{{\bf x}}=(1,0)$, we pick $\tilde{{\bf x}}$ so that 
\[
\tilde{x}_{j} \sim \mathcal{N}(0,\sigma^{2}),
\]
where $\mathcal{N}(0,\sigma^{2})$ denotes a mean zero Gaussian random variable with variance $\sigma^{2}$.  We then generate the reduced dynamical quantities 
\begin{equation}
\mathbb{E}\left[\left.y_{1}(t)\right|\hat{{\bf x}}=(1,0)\right], ~ \mathbb{E}\left[\left.y_{2}(t)\right|\hat{{\bf x}}=(1,0)\right],
\label{redquns}
\end{equation}
by averaging over $N_{ens} = 10^{4}$ randomly generated choices of $\tilde{{\bf x}}$.  Likewise, to generate the results for the MDDMD, we average over an ensemble of $10^{4}$ choices for $K_{0}$ where each component is likewise chosen from $\mathcal{N}(0,\sigma^{2})$.  The unperturbed data forming the data streams $G_{-}$ and $G_{+}$ comes from measuring only $\left(y_{1}(t),y_{2}(t)\right)^{T}$ from one of the random realizations used to generate the quantities in Equation \eqref{redquns}.  We run each simulation for $t_{f}=50$ units of time with a time step of $\delta t =10^{-1}$, which for simplicity we also take to be our sampling rate of the MDDMD, $\delta t_2=10^{-1}$.  

\begin{figure}[H]
\centering
\begin{tabular}{c}
\includegraphics[width=1\textwidth]{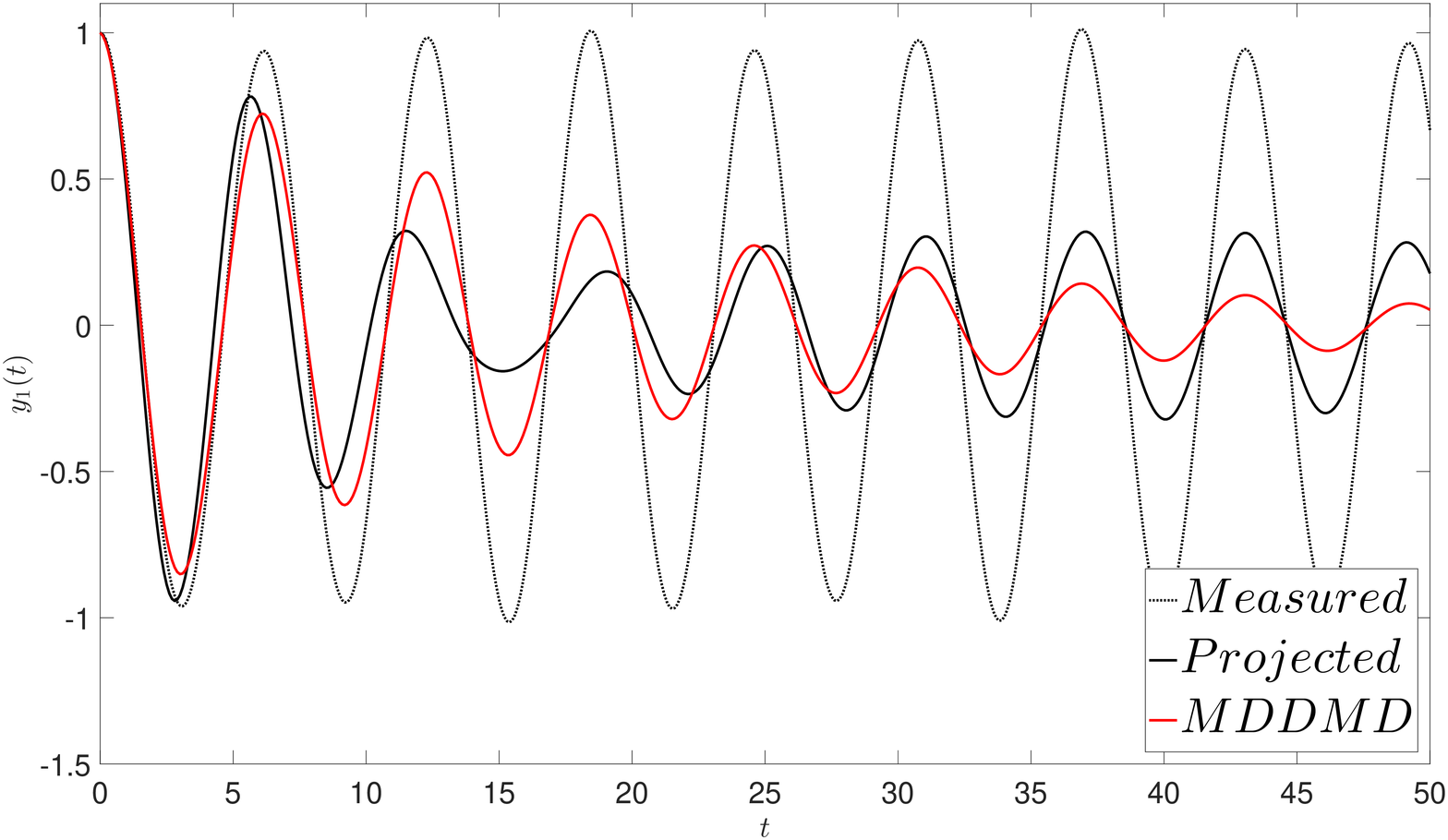}\\
(a) \\
\includegraphics[width=1\textwidth]{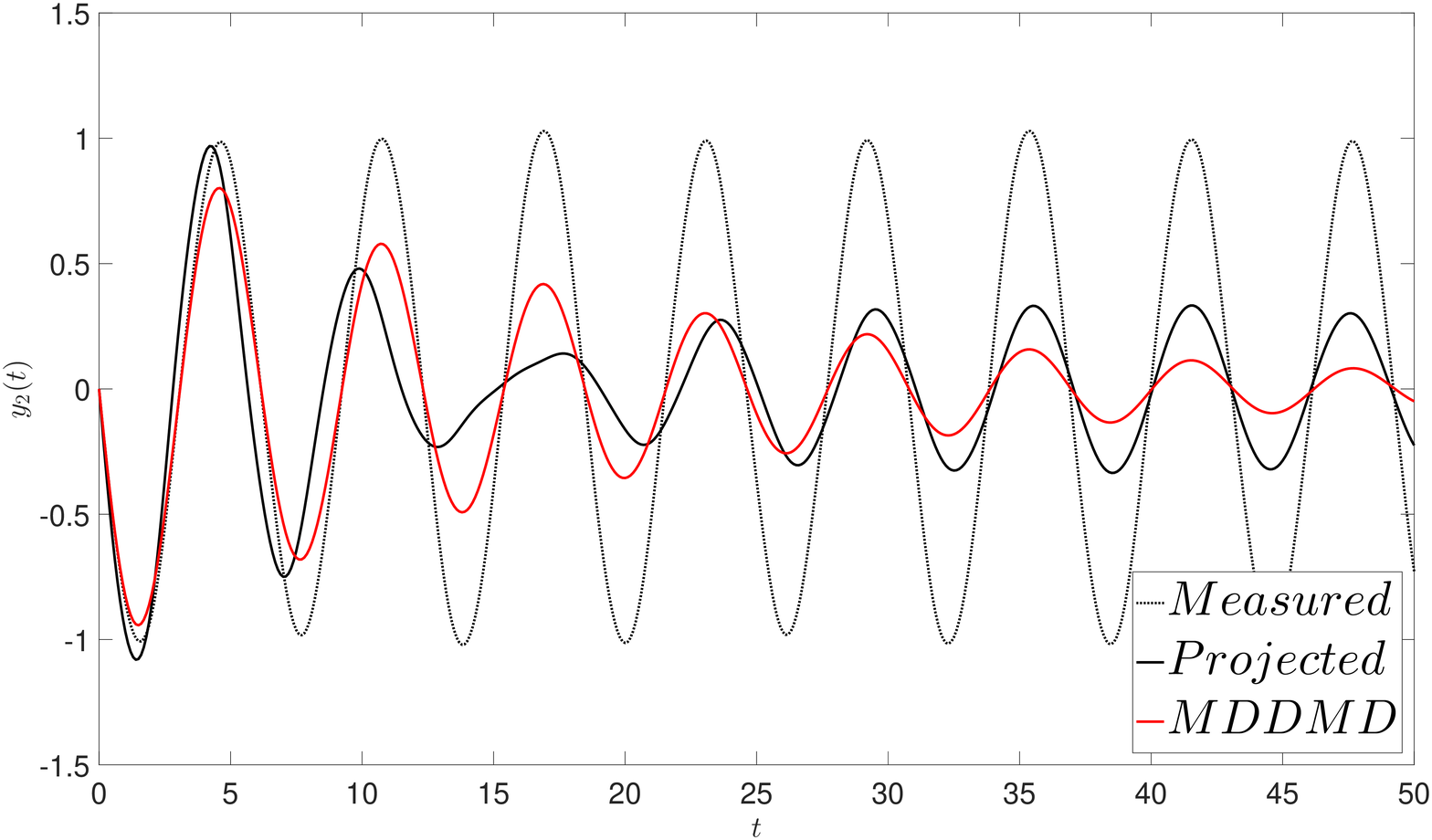}\\
(b)
\end{tabular}
\caption{For $\sigma=.5$, comparison of a randomly generated sample measurement (dotted line), the ensemble reduced dynamics (solid black line), and the MDDMD approximation to the ensemble (solid red line) for the first component (a) and the second component (b).}
\label{fig:sigpt5}
\end{figure}

\begin{figure}[H]
\centering
\begin{tabular}{c}
\includegraphics[width=1\textwidth]{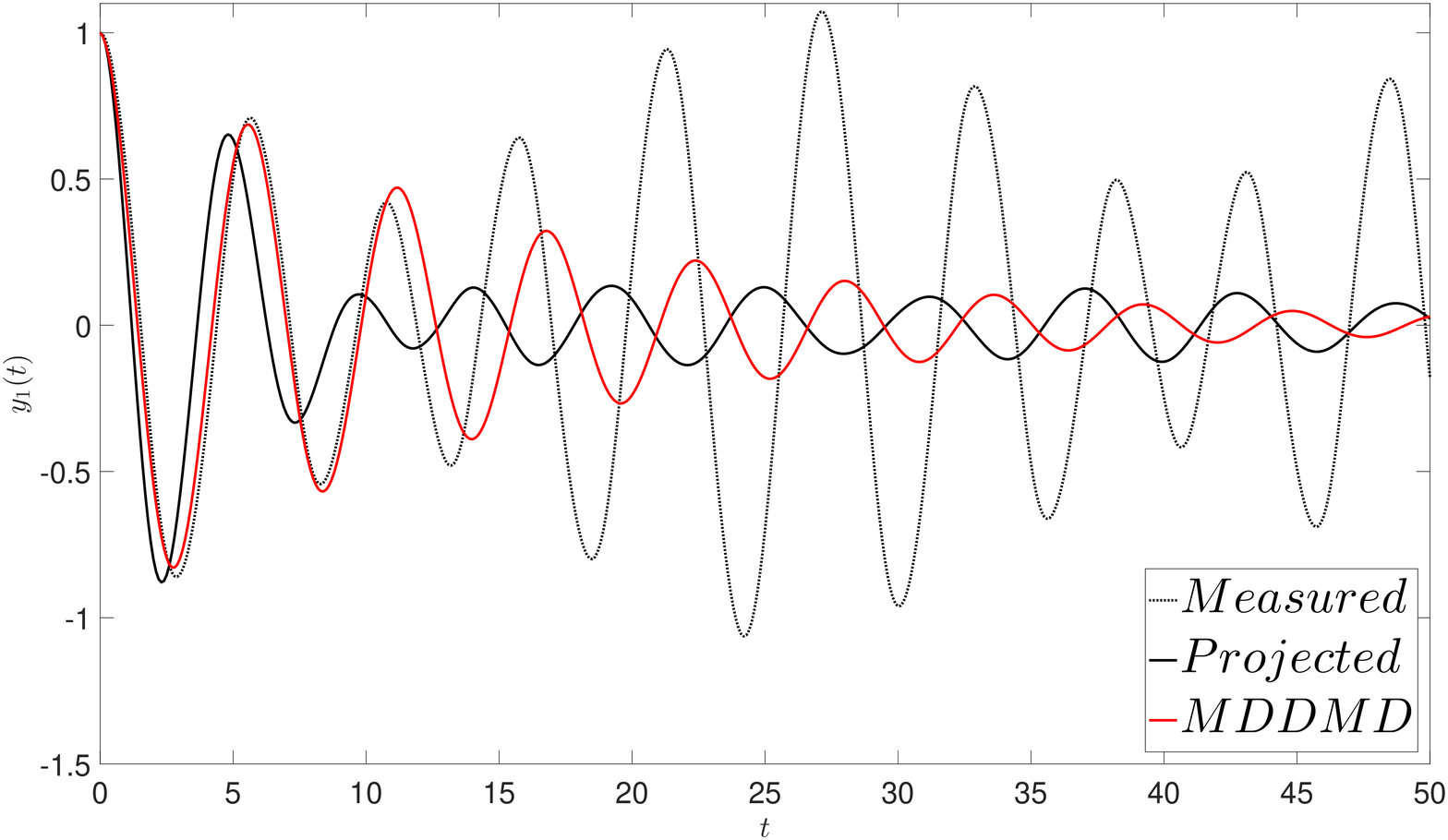}\\
(a) \\
\includegraphics[width=1\textwidth]{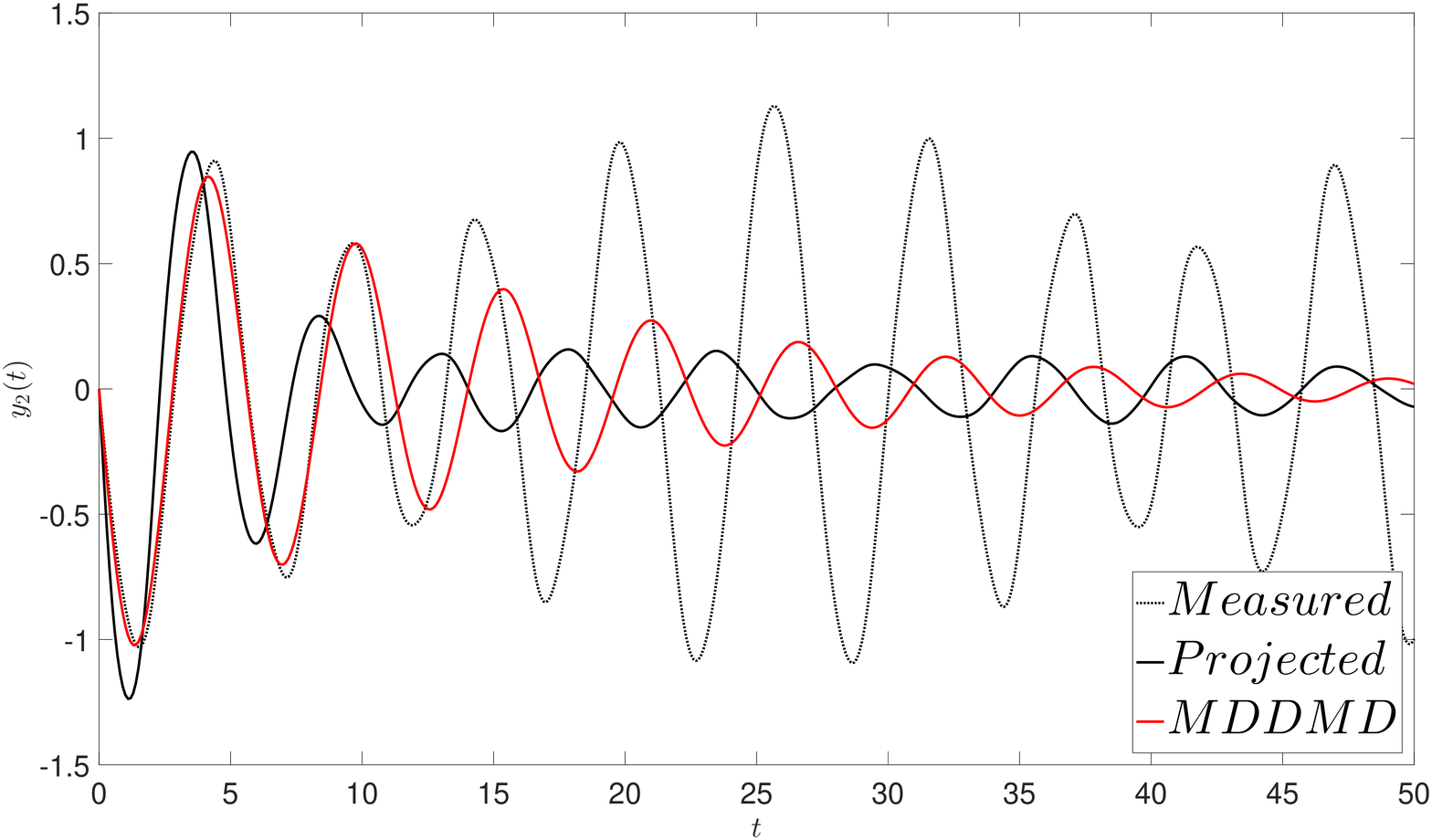}\\
(b)
\end{tabular}
\caption{For $\sigma=1$, comparison of a randomly generated sample measurement (dotted line), the ensemble reduced dynamics (solid black line), and the MDDMD approximation to the ensemble (solid red line) for the first component (a) and the second component (b).}
\label{fig:sig1}
\end{figure}

Setting $\sigma = .5$ generates the results seen in Figure \ref{fig:sigpt5}.  As can be seen, while not exact, the MDDMD is able to reproduce the decaying oscillatory nature of the projected dynamics.  That this is done from a single random measurement which itself does not exhibit the decay properties seen makes the MDDMD results especially compelling.  That said, the typical amplitude of the missing dimensions $(y_{3},y_{4})$ is less than $\sigma$, which makes the coupling to the measured dimensions on the order of $\sigma^{2}$, which is relatively small, and thus the coupling is weak.  Doubling the standard deviation so that $\sigma = 1$, which makes the coupling between measured and unmeasured dimensions significantly stronger generates the results in Figure \ref{fig:sig1}.  While not as accurate, the MDDMD is still able to track much of the reduced dynamics, though the added stochasticity causes the underlying measurement used in the MDDMD to push the approximation out of phase with respect to the true reduced dynamics. 

%%%%%%%%%%%%%%%%%%%%%%%%%%%%%%%%%%%%%%%%%%%%%%%%%%%%%%%%%%%%%%%%%%%%%%%%%%%%%%%%
\subsection{A Slow-Fast Hamiltonian System}
We modify the Hamiltonian system from Section \ref{example_1} so that now
\[
H({\bf y}) = \frac{1}{2}\left( y_{1}^{2}+y_{2}^{2}\right)+\frac{\epsilon}{2}\left(y_{3}^{2}+y_{4}^{2}+ y_{1}^{2}y_{3}^{2}\right), 
\]
which produces the corresponding slow-fast system
\begin{align*}
\dot{y}_1 &= y_2, \\
\dot{y}_2 &= - y_1 \left(1 + \epsilon y_3^2 \right), \\
\dot{y}_3 &= \epsilon y_4, \\
\dot{y}_4 &= -\epsilon y_3 \left(1 + y_1^2 \right).
\end{align*}
Here the unresolved subsystem, $\tilde{{\bf x}}=(x_3,x_4)$, now oscillates at a speed roughly $\epsilon$ times faster than $\hat{{\bf x}}=(x_1,x_2)$. Setting $\epsilon=10$ then simulates a system in which the unresolved variables move an order of magnitude faster, though also for which the coupling between $y_{2}$ and $y_{3}$ is ten times stronger.  

We again use $N_{ens}=10^4$ random initial states for $\tilde{{\bf x}}$, with each component drawn from $\mathcal{N}(0,\sigma^2)$ and subsequently determine the mean dynamics of $\hat{{\bf x}}$ as in Equation \eqref{redquns}.  Given that we are interested in looking at a system in which the missing dimensions are fast, but not necessarily of large magnitude, we require that $\sigma \leq 1/\sqrt{\epsilon}$.   To generate the exactly averaged process, we must decrease the time step in our numerical simulations to $\delta t=10^{-2}$ in order to accurately compute the fast-moving subsystem. However, the measurement sample rate, or MDDMD sampling rate, is left unchanged so that $\delta t_{2}=10^{-1}$.

\begin{figure}[H]
\centering
\begin{tabular}{c}
\includegraphics[width=1\textwidth]{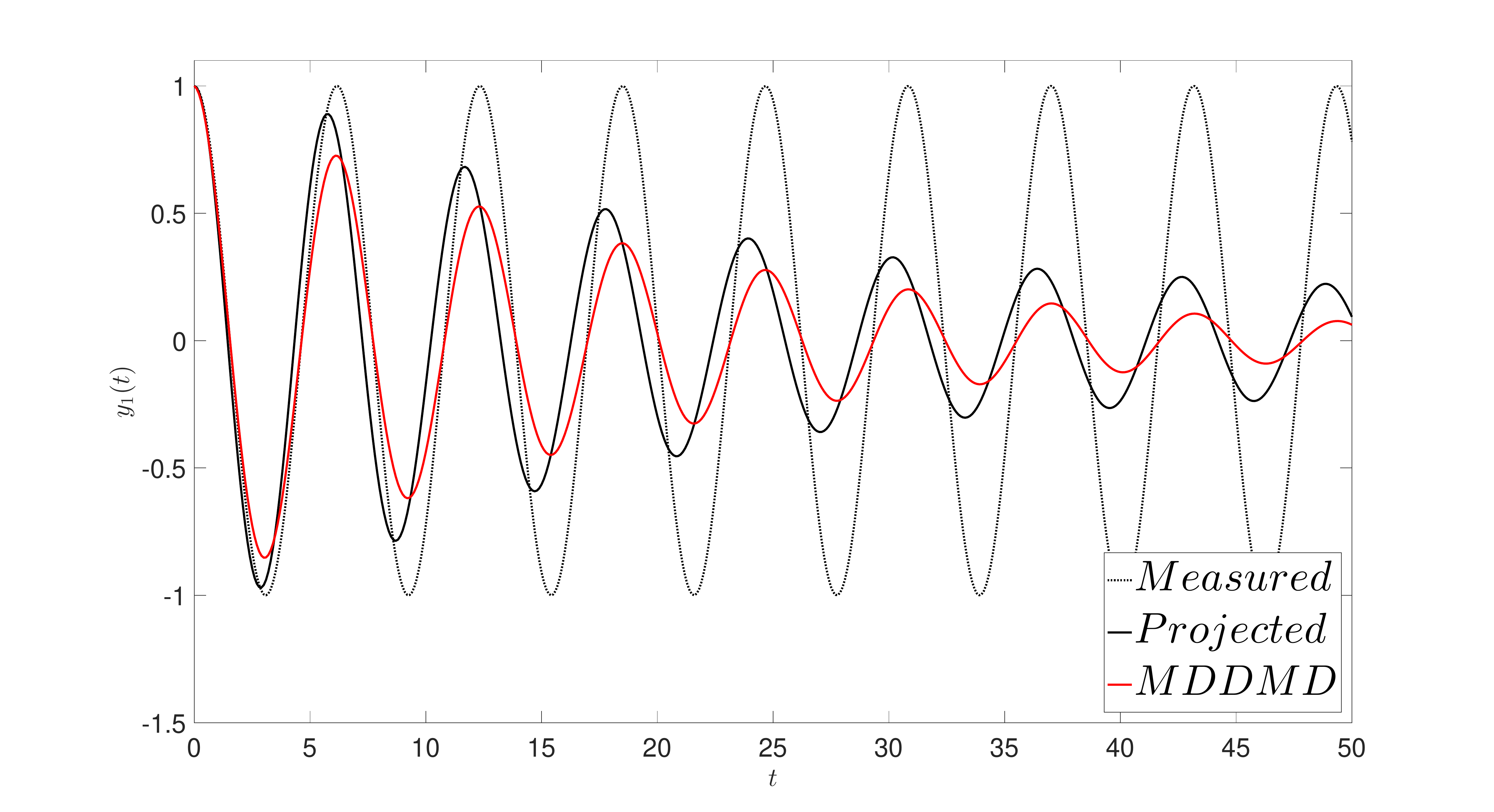}\\
(a) \\
\includegraphics[width=1\textwidth]{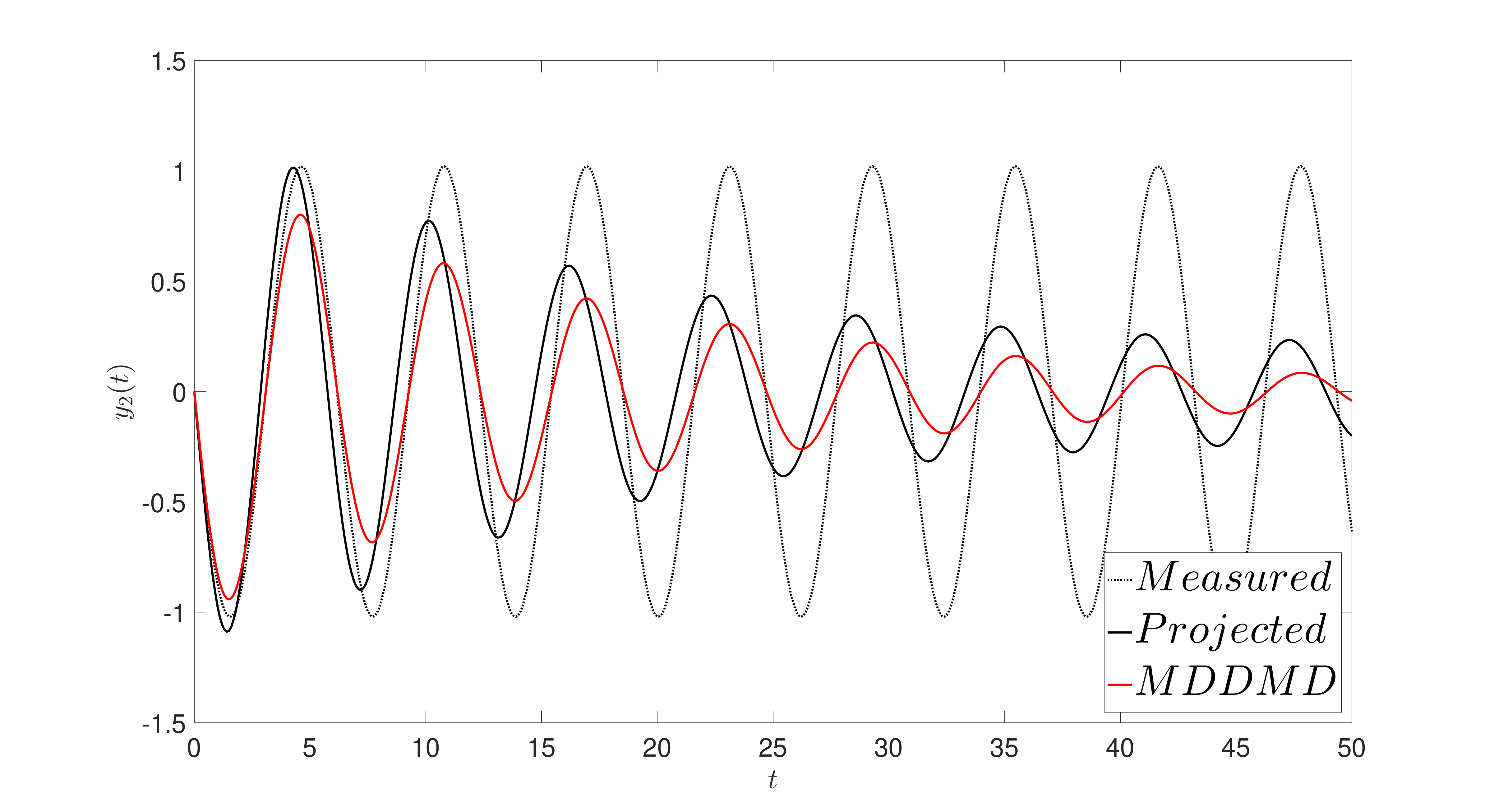}\\
(b)
\end{tabular}
\caption{For $\sigma=0.5/\sqrt{\epsilon}$ and $\epsilon=10$, comparison of a randomly generated sample measurement (dotted line), the ensemble reduced dynamics (solid black line), and the MDDMD approximation to the ensemble (solid red line) for the first component (a) and the second component (b).}
\label{fig:slowfast-sigpt5}
\end{figure}

\begin{figure}[H]
\centering
\begin{tabular}{c}
\includegraphics[width=1\textwidth]{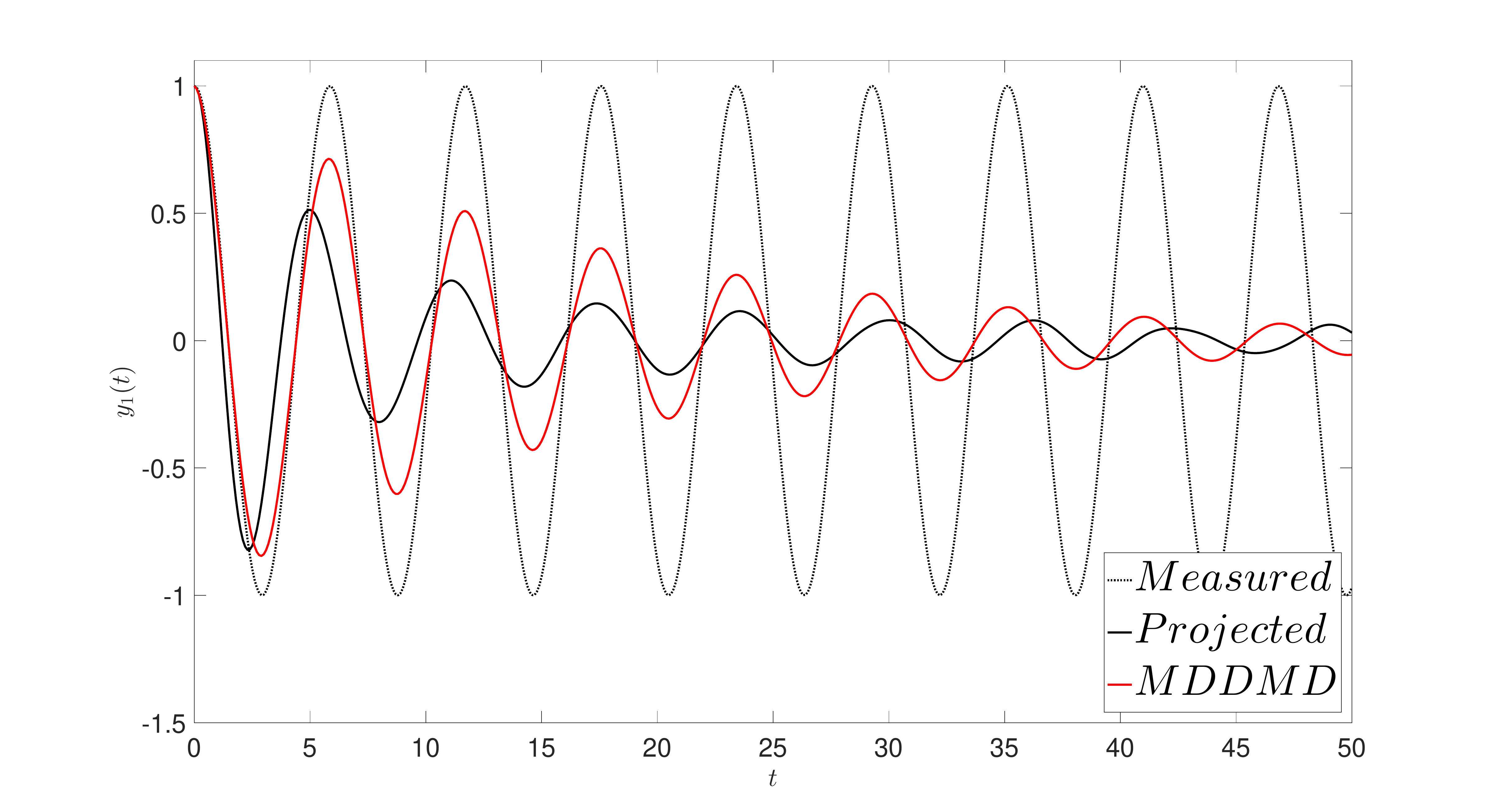}\\
(a) \\
\includegraphics[width=1\textwidth]{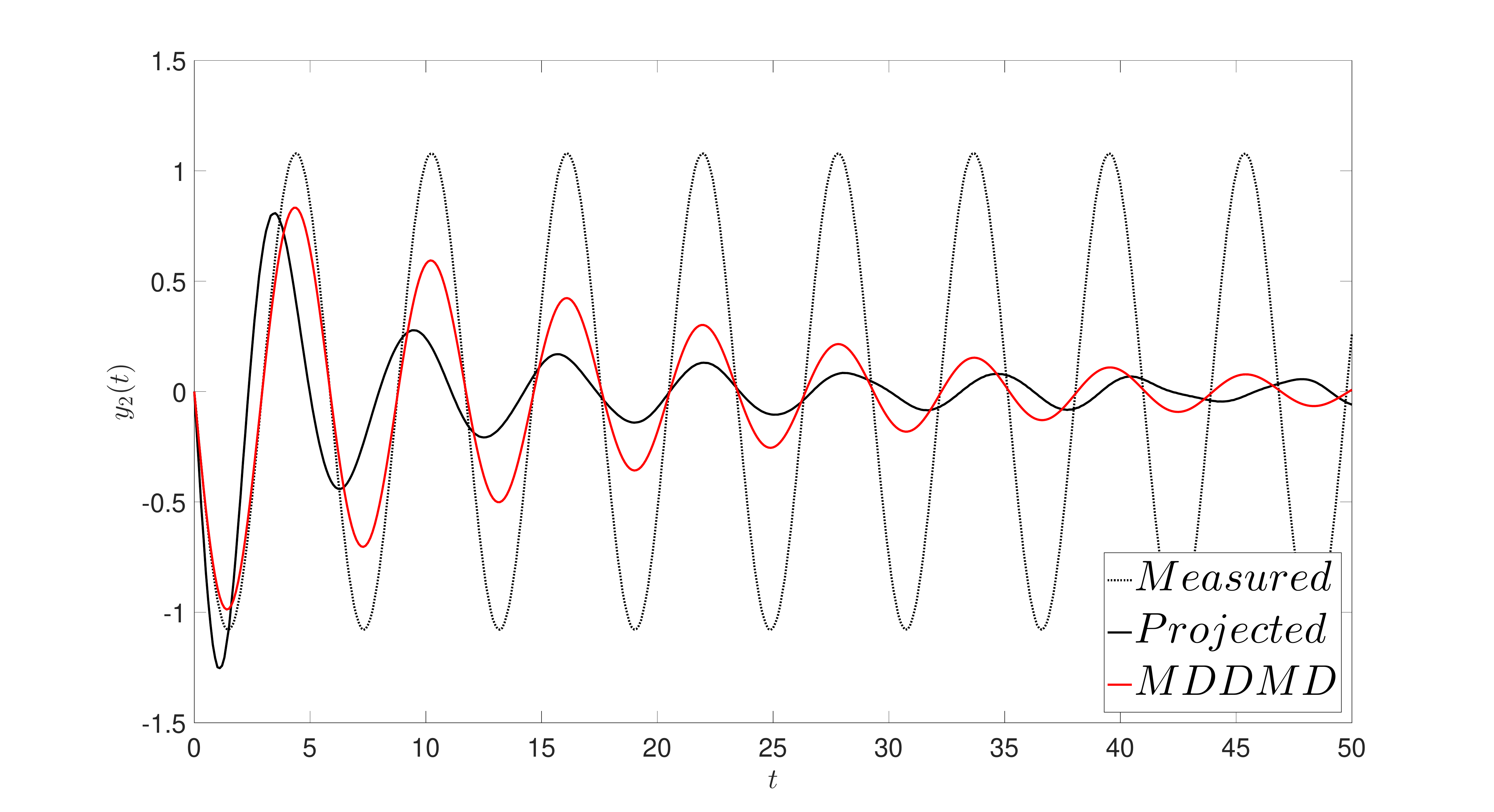}\\
(b)
\end{tabular}
\caption{For $\sigma=1/\sqrt{\epsilon}$ and $\epsilon=10$, comparison of a randomly generated sample measurement (dotted line), the ensemble reduced dynamics (solid black line), and the MDDMD approximation to the ensemble (solid red line) for the first component (a) and the second component (b).}
\label{fig:slowfast-sig1}
\end{figure}

For $\sigma=0.5/\sqrt{\epsilon}$, Figure \ref{fig:slowfast-sigpt5} demonstrates how the MDDMD approximation can still readily determine the rate of decay of the projected dynamics even in the presence of fast-moving unresolved variables. Figure \ref{fig:slowfast-sig1} demonstrates similar results for the $\sigma=1/\sqrt{\epsilon}$ case. However, it is prudent to note that with $\epsilon=10$, the fast unresolved subsystem $\tilde{{\bf x}}$ begins to average out quicker over the Monte Carlo simulations. Thus, the projected dynamics are mostly dominated by their decay towards equilibrium and less sensitive to the fast subsystem,  yet MDDMD is still able to determine the appropriate rate of decay given only a single non-decaying observation.
%\begin{figure}[H]
%\centering
%\begin{tabular}{c}
%\includegraphics[width=1\textwidth]{y12_slow}\\
%(a) \\
%\includegraphics[width=1\textwidth]{y34_fast}\\
%(b)
%\end{tabular}
%\caption{For $\epsilon=10$, comparison of a single realization of the measured slow subsystem (a) and the randomly generated fast subsystem (b).}
%\label{fig:slowfast}
%\end{figure}

\begin{figure}[H]
\centering
\begin{tabular}{c}
\includegraphics[width=1\textwidth]{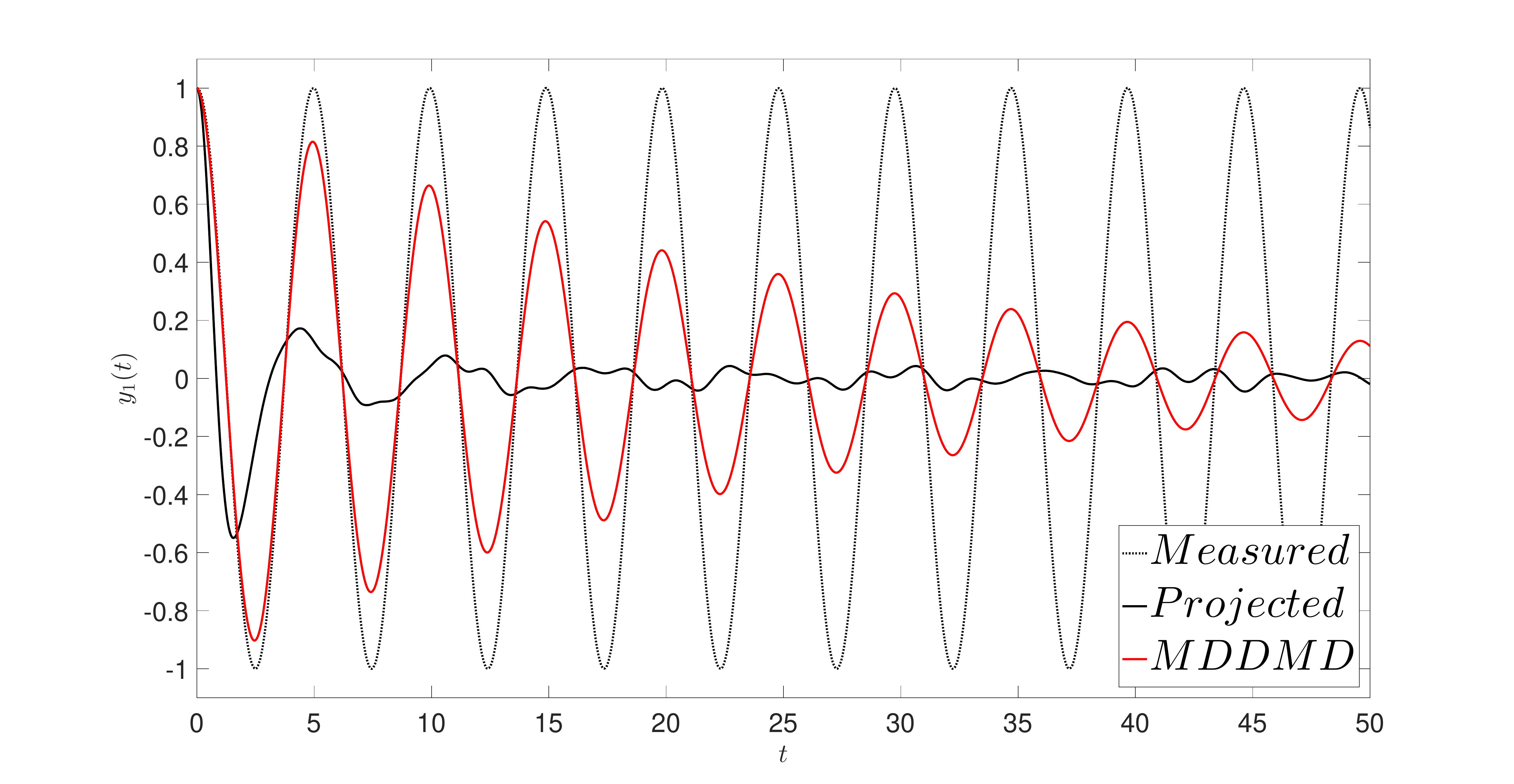}\\
(a) \\
\includegraphics[width=1\textwidth]{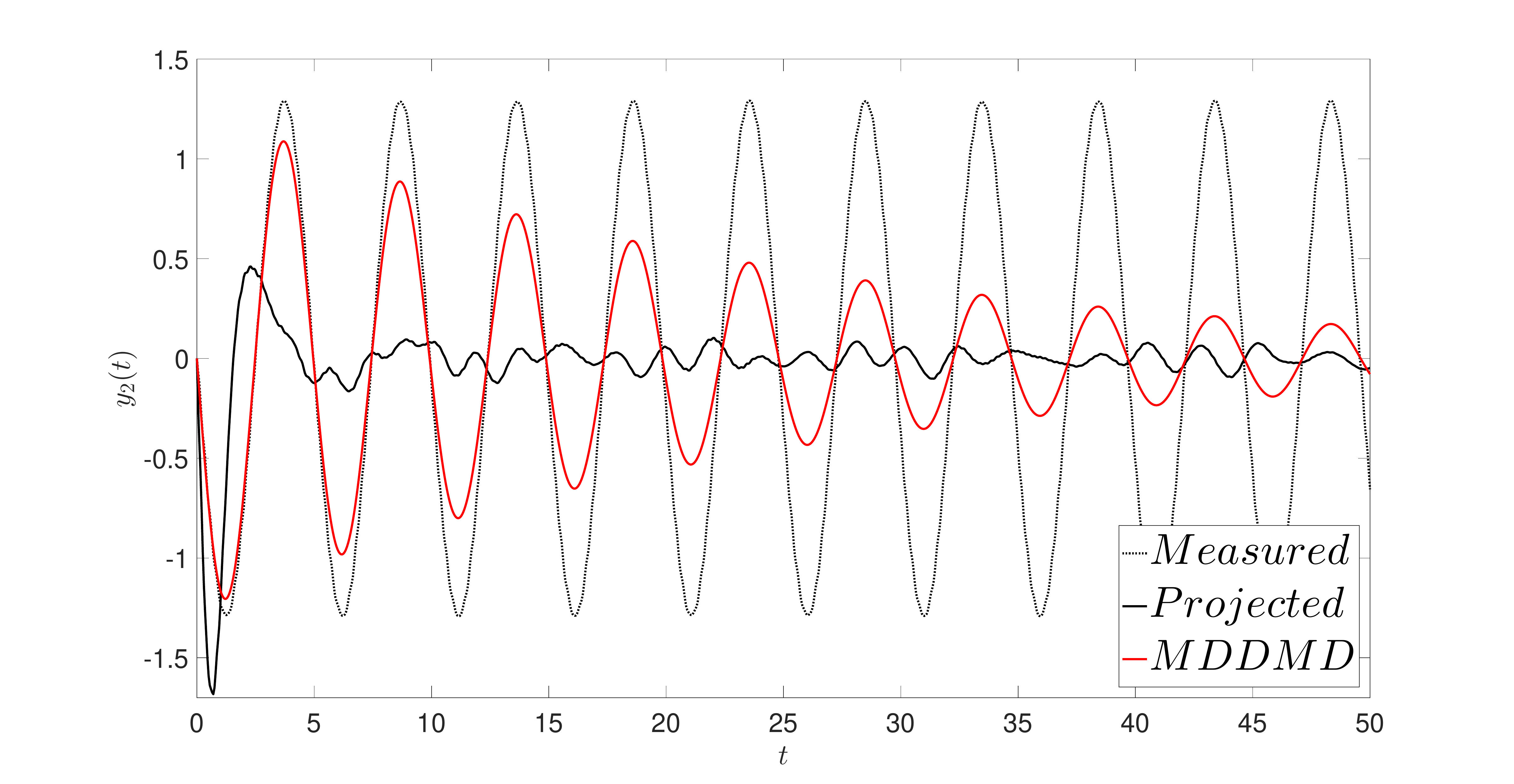}\\
(b)
\end{tabular}
\caption{For $\sigma=2/\sqrt{\epsilon}$ and $\epsilon=10$, comparison of a randomly generated sample measurement (dotted line), the ensemble reduced dynamics (solid black line), and the MDDMD approximation to the ensemble (solid red line) for the first component (a) and the second component (b).}
\label{fig:slowfast-sig2}
\end{figure}

It is interesting to examine a case outside of the asymptotic regime in which we might expect our method to work.   Thus, we let $\sigma=2/\sqrt{\epsilon}$, which now allows for a far greater chance of initial choices in $(y_{3},y_{4})$ which lead to significantly stronger couplings between the measured and unmeasured dimensions.  As seen in Figure \ref{fig:slowfast-sig2}, the combination of faster time scale and larger amplitude makes it difficult for the MDDDM method to faithfully capture even the rate of decay the true average.  Moreover, the faster and stronger oscillations in the missing dimensions add corresponding greater degrees of oscillation in the true mean.  Addressing these shortcomings of the MDDMD method is a subject for future research.  

%%%%%%%%%%%%%%%%%%%%%%%%%%%%%%%%%%%%%%%%%%%%%%%%%%%%%%%%%%%%%%%%%%%%%%%%%%%%%%%%
% Conclusion and Future Directions
%%%%%%%%%%%%%%%%%%%%%%%%%%%%%%%%%%%%%%%%%%%%%%%%%%%%%%%%%%%%%%%%%%%%%%%%%%%%%%%%
\section{Conclusion and Future Directions}
We have derived in this paper an extension to the original DMD algorithm, the MDDMD algorithm, which facilitates accounting for the effects of missing information.  This is done by introducing a memory kernel in the vein of what has been done in the statistical mechanics literature.  This allowed us then to develop a method which produces reasonable approximations to averaged dynamics from one incomplete measurement.  This likewise permits us to generate a time-predictive model accounting for the impact of missing information in a dynamical system.  While the use of our method on a toy problem was studied in this paper, given the relative success, it is a future project to use this in more complex models or on real-world data.  Likewise, while a reasonably self-consistent approximation scheme was used to derive the MDDMD approach, it is an especially interesting question to explore other approximations which facilitate better resolving the impact of orthogonal dynamics, such as in \cite{stinis}.  This is also a subject of future research.  

%%%%%%%%%%%%%%%%%%%%%%%%%%%%%%%%%%%%%%%%%%%%%%%%%%%%%%%%%%%%%%%%%%%%%%%%%%%%%%%%
% Appendix
%%%%%%%%%%%%%%%%%%%%%%%%%%%%%%%%%%%%%%%%%%%%%%%%%%%%%%%%%%%%%%%%%%%%%%%%%%%%%%%%
\section{Appendix}
%%%%%%%%%%%%%%%%%%%%%%%%%%%%%%%%%%%%%%%%%%%%%%%%%%%%%%%%%%%%%%%%%%%%%%%%%%%%%%%%
\subsection{Deriving the MZD}
We briefly collect some notes on the derivation of the Mori-Zwanzig Decomposition (MZD) given by Equations \eqref{mzeq1} and \eqref{mzeq2} of the projected dynamical system.  This starts by looking at the affiliated evolution equation
\[
w_{t} = \mathbb{Q}\mathcal{L}w, ~ w\left({\bf x},0\right) = w_{0}\left({\bf x}, \right)
\]
using 
\[
\pd_{t}\left(e^{-\mathcal{L}t}w \right) = -e^{-\mathcal{L}t}\mathbb{P}\mathcal{L}w,
\]
we see that 
\[
w({\bf x},t)  = e^{\mathcal{L}t}w_{0} - \int_{0}^{t}e^{\mathcal{L}(t-s)}\mathbb{P}\mathcal{L}w({\bf x},s)~ds,
\]
so that after rearrangement we derive Dyson's formula
\[
e^{\mathcal{L}t} = e^{\mathbb{Q}\mathcal{L}t} + \int_{0}^{t}e^{\mathcal{L}(t-s)}\mathbb{P}\mathcal{L}e^{\mathbb{Q}\mathcal{L}s}~ds.
\]
Using this on $\mathbb{Q}\mathcal{L}g$ provides the identity
\[
e^{\mathcal{L}t}\mathbb{Q}\mathcal{L}g = e^{\mathbb{Q}\mathcal{L}t}\mathbb{Q}\mathcal{L}g + \int_{0}^{t}e^{\mathcal{L}(t-s)}\mathbb{P}\mathcal{L}e^{\mathbb{Q}\mathcal{L}s}\mathbb{Q}\mathcal{L}g ~ds,
\]
which then lets us decompose $\pd_{t}\left(e^{\mathcal{L}t}g \right) = e^{\mathcal{L}t}\mathcal{L}g$ for general observable $g({\bf x})$ as 
\begin{align*}
e^{\mathcal{L}t}\mathcal{L}g = & e^{\mathcal{L}t}\mathbb{P}\mathcal{L}g + e^{\mathcal{L}t}\mathbb{Q}\mathcal{L}g,\\
= & e^{\mathcal{L}t}\mathbb{P}\mathcal{L}g +  e^{\mathbb{Q}\mathcal{L}t}\mathbb{Q}\mathcal{L}g + \int_{0}^{t}e^{\mathcal{L}(t-s)}\mathbb{P}\mathcal{L}e^{\mathbb{Q}\mathcal{L}s}\mathbb{Q}\mathcal{L}g ~ds.
\end{align*}
Given the noise $F$ to be 
\[
F(t;g) =  e^{\mathbb{Q}\mathcal{L}t}\mathbb{Q}\mathcal{L}g,
\]
and given that for all reasonable choices of functions $\tilde{f}$ we have the dynamical orthogonality condition 
\[
\mathbb{P}e^{\mathbb{Q}\mathcal{L}t}\mathbb{Q}\tilde{f}= 0,
\]
we finally have the reduced equation
\[
\mathbb{P}e^{\mathcal{L}t}\mathcal{L}g = \mathbb{P} e^{\mathcal{L}t}\mathbb{P}\mathcal{L}g + \int_{0}^{t}\mathbb{P}e^{\mathcal{L}(t-s)}K(s;g) ~ds,
\]
or we find Equation \eqref{mzeq1}
\[
\pd_{t}\left(\mathbb{P}e^{\mathcal{L}t} g \right) = \mathbb{P} e^{\mathcal{L}t}\mathbb{P}\mathcal{L}g + \int_{0}^{t}\mathbb{P}e^{\mathcal{L}(t-s)}K(s;g) ~ds.
\]
where the `memory kernel' $K(s;g)$ is given by
\[
K(s;g) = \mathbb{P}\mathcal{L}F(s;g).
\]
We likewise have the complementary orthogonal noise equation
\[
\pd_{t}F = \mathbb{Q}\mathcal{L}F, ~ F(0;g) = \mathbb{Q}\mathcal{L}g,
\]
so that by using Dyson's formula, we have that
\[
F(t;g) = e^{\mathcal{L}t}F(0;g) - \int_{0}^{t}e^{\mathcal{L}(t-s)}K(s;g) ~ ds,
\]
and so, by multiplying by $\mathbb{P}\mathcal{L}$, we then Equation \eqref{mzeq2}
\[
K(t;g) + \int_{0}^{t}\mathbb{P}e^{\mathcal{L}(t-s)}\mathcal{L}K(s;g) ~ ds= \mathbb{P}e^{\mathcal{L}t}\tilde{F}(0;g).
\]
where $\tilde{F}(0;g) = \mathcal{L}F(0;g)$.  

%%%%%%%%%%%%%%%%%%%%%%%%%%%%%%%%%%%%%%%%%%%%%%%%%%%%%%%%%%%%%%%%%%%%%%%%%%%%%%%%
\subsection{Finding $DM(\tilde{\mathcal{K}}_{a,0};K_0)$}
In order to find $DM(\tilde{\mathcal{K}}_{a,0};K_0)$, we must compute
\begin{multline*}
\mbox{tr}\left(R^{T}\big(\tilde{\mathcal{K}}_{a,0}\big)D\tilde{M}(\tilde{K}_{a,0};K_{0})W\right) = \\
\lim_{\epsilon\rightarrow 0}\mbox{tr}\left( R^{T}\big(\tilde{\mathcal{K}}_{a,0}\big)\frac{\left(\tilde{M}\big(\tilde{\mathcal{K}}_{a,0}+\epsilon W; K_0\big) - \tilde{M}\big(\tilde{\mathcal{K}}_{a,0};K_0\big)\right)}{\epsilon} \right).
\end{multline*}
First, let 
\[
\tilde{A} = \tilde{\mathcal{K}}_{a,0} - I.
\]
Then setting $S(\tilde{A}) = \left(1 - \frac{\delta t}{2}\right)I + \frac{1}{2}\tilde{A}$, we rewrite the $\tilde{M}(\tilde{A})$ matrix from Section 2 as
\[
\tilde{M}(\tilde{A};K_0) = S(\tilde{A}) \hat{M}(\tilde{A}),
\]
where
\[
\hat{M}(\tilde{A}) = \{0 f_1(\tilde{A})K_0 \cdots f_n(\tilde{A})K_0\}.
\]
Thus, we find immediately that 
\[
\mbox{tr}\left(R^{T}D\tilde{M}W\right) = \frac{1}{2}\mbox{tr}\left(\hat{M}R^{T}W\right) + \lim_{\epsilon\rightarrow 0} J_{\hat{M}}(\tilde{A},W,\epsilon),
\]
where
\begin{equation}
J_{\hat{M}}(\tilde{A},W,\epsilon) = \mbox{tr}\left( R^{T}(\tilde{K}_{a,0})S(\tilde{A})\left( \frac{\hat{M}\left(\tilde{A}+\epsilon W\right) - \hat{M}(\tilde{A})}{\epsilon} \right) \right). \label{dmw_eqn}
\end{equation}
Choosing now the particular perturbation matrix, $W^{(jk)}$, such that
\[ W^{(jk)}_{mn} =
  \begin{cases}
    1 &\quad \text{if } m=j, \text{ and } k=n \\
    0 &\quad \text{if } m\ne j, \text{ or } k \ne n,
  \end{cases}
\]
we then have
\[
\mbox{tr}\left(R^{T}(\tilde{K}_{a,0})S(\tilde{A})D\hat{M}W^{(jk)}\right) = (R^{T}(\tilde{K}_{a,0})S(\tilde{A})D\hat{M})_{kj}.
\]
Now assuming there is a suitable diagonalization of the matrix $\tilde{A}$,
\[
\tilde{A} = V_0\Lambda_0 V_0^{-1},
\]
$\hat{M}(\tilde{A})$ then becomes
\[
\hat{M}(\tilde{A}) = V_0 \Big\{0 \  f_1(\Lambda_0)\tilde{K}_0 \ \dots \ f_n(\Lambda_0)\tilde{K}_0 \Big\},
\]
where
\[
\tilde{K}_0 = V_0^{-1} K_0.
\]
For the matrix $\tilde{M}(\tilde{A} + \epsilon W^{(jk)})$ we need to first find a perturbative approximation for the diagonalization of $\tilde{A} + \epsilon W^{(jk)}$. To wit, we see that
\begin{align*}
\tilde{A} + \epsilon W^{(jk)} &= V \left(\Lambda_0 + \epsilon \tilde{W}^{(jk)}\right)V^{-1}, \\
\tilde{W}^{(jk)} &= V^{-1}W^{(jk)}V.
\end{align*}
Solving the eigenvalue problem,
\[
\left(\Lambda_0 + \epsilon\tilde{W}^{(jk)}\right)(\phi_{0,l}+\epsilon\phi_{1,l}+\dots)=(\lambda_{0,l}+\epsilon\lambda_{1,l}+\dots)(\phi_{0,l}+\epsilon\phi_{1,l} +\dots),
\]
we then have
\[
\Lambda_0\phi_{0,l}=\lambda_{0,l} \phi_{0,l}, 
\]
and to next order
\[
(\Lambda_0 - \lambda_{0,l})\phi_{1,l} = \left(\lambda_1 - \tilde{W}^{(jk)} \right) \phi_{0,l}.
\]
For the sake of simplicity, we suppose that the eigenvalues of $\tilde{A}$ are simple, so that $\phi_{0,l}=\hat{e}_l$, i.e. the $l^{th}$ canonical basis vector. We then find that
\[
\lambda_{1,l} = \left< \hat{e}_l, \tilde{W}^{(jk)}\hat{e}_l \right>,
\]
and likewise, that
\[
\phi_{1,l}=(\Lambda_0 - \lambda_{0,l})^{-P} \left(\lambda_{1,l} - \tilde{W}^{(jk)} \right) \hat{e}_l,
\]
where $X^{-P}$ is the Moore-Penrose pseudoinverse for the matrix $X$. Thus, we find that to next order
\[
\tilde{A} + \epsilon W^{(jk)} = V_0 (I + \epsilon \Phi_1)(\Lambda_0 + \epsilon \Lambda_1) (I + \epsilon \Phi_1)^{-1} V_0^{-1}.
\]
Hence,
\begin{multline*}
\hat{M}\left(\tilde{A} + \epsilon W^{(jk)}\right) = \\
V_0 (I +\epsilon\Phi_1) \Big\{ 0 \ f_1(\Lambda_0 + \epsilon\Lambda_1)(I+\epsilon\Phi_1)^{-1}\tilde{K}_0 \ \cdots \ f_n(\Lambda_0 + \epsilon\Lambda_1)(I+\epsilon\Phi_1)^{-1}\tilde{K}_0 \Big\}.
\end{multline*}
Using
\[
(I+\epsilon \Phi_1)^{-1} = I - \epsilon\Phi_1 + \mathcal{O}(\epsilon^2),
\]
we then get
\begin{align*}
\hat{M}\left( \tilde{A} + \epsilon W^{(jk)} \right) &= V_0 \Big\{0 \ f_1(\Lambda_0 + \epsilon\Lambda_1)\tilde{K}_0 \ \cdots \ f_n(\Lambda_0 + \epsilon\Lambda_1)\tilde{K}_0 \Big\} \\
&\ + \epsilon V_0 \Big\{ 0 ~[\Phi_1,f_1(\Lambda_0)]\tilde{K}_0 \ \cdots \ [\Phi_1,f_n(\Lambda_0)]\tilde{K}_0 \Big\}, 
\end{align*}
where $[\cdot, \cdot]$ is the commutator. Hence, we get the final approximation
\begin{flalign*}
&\frac{1}{\epsilon} \left( \hat{M}\left( \tilde{A} + \epsilon W^{(jk)} \right) - \hat{M}\left(\tilde{A}\right)\right) \sim \\
&\quad\quad \frac{1}{\epsilon}V_0 \Big\{ 0 \ (f_1(\Lambda_0 + \epsilon\Lambda_1) - f_1(\Lambda_0))\tilde{K}_0 \ \cdots \ (f_n(\Lambda_0 + \epsilon\Lambda_1)-f_n(\Lambda_0))\tilde{K}_0 \Big\} \\
&\quad\quad + V_0 \Big\{ 0 ~[\Phi_1,f_1(\Lambda_0)]\tilde{K}_0 \ \cdots \ [\Phi_1,f_n(\Lambda_0)]\tilde{K}_0 \Big\}.
\end{flalign*}
Updating this matrix and substituting back into equation (\ref{dmw_eqn}) for each $W^{(jk)}$ perturbation then completes the algorithm for finding $DM(\tilde{A};K_0)$. Hence, by exploiting the diagonalization of the matrix $\tilde{A}$, we are able to go from having to solve $N^2$ eigenvalue problems to only computing $N^2$ traces, a much simpler operation.

\bibliography{dmd_wavelet}
\bibliographystyle{unsrt}
\end{document}